\documentclass[11pt]{article}

\newtheorem{lem}{Lemma}
\newtheorem{theo}{Theorem}

\newtheorem{prop}{Proposition}

\newenvironment{pf}{\ \\ {\bf Proof: }}{\hfill\mbox{$\diamond$}\medskip}
\usepackage{amsmath}
\usepackage{amssymb}

\begin{document}
\author{Vlad Bally\thanks{Laboratoire d'Analyse et de
Math\'ematiques Appliqu\'ees, UMR $8050$,
Universit\'e Paris-Est Marne-la-Vall\'ee, 5  Bld Descartes,
Champs-sur-marne, 77454 Marne-la-Vall\'ee Cedex 2, France.} \footnote{  
Acknowledgement : Part of this work has been done during a visit of
the first author to the Institute Mittag-Leffler (Djursholm, Sweden)}
\and  Emmanuelle Cl\'ement$^*$}
\title{Integration by parts formula with respect to jump times for stochastic differential equations}
\date{April, 7 2010}
\maketitle
\begin{abstract}
We establish an integration by parts formula based on jump times in an abstract framework in order to study the regularity 
of the law for processes solution of stochastic differential equations with jumps.



\end{abstract}

\section{Introduction}

We consider the one dimensional equation 
\begin{equation}
X_{t}=x+\int_{0}^{t}c(u,a,X_{u-})dN(u,a)+\int_{0}^{t}g(u,X_{u})du  \label{i1}
\end{equation}
where $N$ is a Poisson point measure of intensity measure $\mu $ on some
abstract measurable space $E.$ We assume that $c$ and $g$ are infinitely
differentiable with respect to $t$ and $x$, have bounded derivatives of any
order and have linear growth with respect to $x.$ Moreover we assume that
the derivatives of $c$ are bounded by a function $\overline{c}$ such that $
\int_E \overline{c}(a)d\mu(a) <\infty .$ Under these hypotheses the equation has a
unique solution  and the stochastic integral with respect to the Poisson point measure is a Stieltjes
integral.

Our aim is to give sufficient conditions in order to prove that the law of $X_{t}$ is
absolutely continuous with respect to the Lebesgue measure and has a smooth
density. If $E=\mathbb{R}^{m}$ and if the measure $\mu$ admits a smooth density $h$ then
one may develop a Malliavin calculus based on the amplitudes of the jumps in
order to solve this problem. This has been done first in \cite{Bi} and then in
\cite{BGJ}. But if $\mu $ is a singular measure this approach fails and one has
to use the noise given by the jump times of the Poisson point measure in
order to settle a differential calculus analogous to the Malliavin calculus. This is a much more
delicate problem and several approaches have been proposed. A first step is
to prove that the law of $X_{t}$ is absolutely continuous with respect to
the Lebesgue measure, without taking care of the regularity. A first result
in this sense was obtained by Carlen and Pardoux in \cite{CP} and was followed
by a lot of other papers (see \cite{D}, \cite{ET}, \cite{K1}, \cite{NS}).\ The second step is to
obtain the regularity of the density. Recently two results of this type have
been obtained by Ishikawa and Kunita in \cite{IK} and by Kulik in \cite{K2}. In both
cases one deals with an equation of the form 
\begin{equation}
dX_{t}=g(t,X_{t})dt+f(t,X_{t-})dU_{t}  \label{i2}
\end{equation}
where $U$ is a L\'{e}vy process. The above equation is multi-dimensional
(let us mention that the method presented in our paper may be used in the
multi-dimensional case as well, but then some technical problems related to
the control of the Malliavin covariance matrix have to be solved - and for
simplicity we preferred to leave out this kind of difficulties in this
paper). Ishikawa and Kunita in \cite{IK} used the finite difference approach
given by J.\ Picard in \cite{P} in order to obtain sufficient conditions for the
regularity of the density of the solution of an equation of type (\ref{i1})
(in a somehow more particular form, closed to linear equations). The result
in that paper produces a large class of examples in which we get a smooth
density even for an intensity measure which is singular with respect to the
Lebesgue measure. The second approach is due to Kulik \cite{K2}. He settled a
Malliavin type calculus based on perturbations of the time structure in
order to give sufficient conditions for the smoothness of the density. In
his paper the coefficient $f$ is equal to one so the non degeneracy comes
from the drift term $g$ only. As before, he obtains the regularity of the
density even if the intensity measure $\mu $ is singular. He also proves
that under some appropriate conditions, the density is not smooth for a
small $t$ so that one has to wait before the regularization effect of the
noise produces a regular density.

The result proved in our paper is the following. We consider the function 
\begin{equation*}
\alpha (t,a,x)=g(x)-g(x+c(t,a,x))+(g\partial _{x}c+\partial _{t}c)(t,a,x).
\end{equation*}
Except the regularity and boundedness conditions on $g$ and $c$ we consider
the following non degeneracy assumption. There exists a measurable function $
\underline{\alpha }$ such that $\left\vert \alpha (t,a,x)\right\vert \geq 
\underline{\alpha }(a)>0$ for every $(t,a,x)\in \mathbb{R}_{+}\times E\times \mathbb{R}.$ We
assume that there exists a sequence of subsets $E_{n}\uparrow E$ such that $
\mu (E_{n})<\infty $ and 
\begin{equation*}
\underline{\lim }_{n\rightarrow \infty }\frac{1}{\mu (E_{n})}\ln
(\int_{E_{n}}\frac{1}{\underline{\alpha }(a)}d\mu (a))=\theta <\infty .
\end{equation*}
If $\theta =0$ then, for every $t>0,$ the law of $X_{t}$ has a $C^{\infty }$
density with respect to the Lebesgue measure. Suppose now that $\theta >0$
and let $q\in \mathbb{N}.$ Then, for $t>16\theta (q+2)(q+1)^{2}$ the law of $X_{t}$
has a density of class $\mathcal{C}^{q}.$ Notice that for small $t$ we are not able to
prove that a density exists and we have to wait for a sufficiently large $t$
in order to obtain a regularization effect. 

In the paper of Kulik \cite{K2} one takes $c(t,a,x)=a$ so $\alpha
(t,a,x)=g(x)-g(x+c(t,a,x)).$ Then the non degeneracy condition concerns just
the drift coefficient $g.$ And in the paper of Ishikawa and Kunita the basic
example (which corresponds to the geometric L\'{e}vy process) is $
c(t,a,x)=xa(e^{a}-1)$ and $g$ constant. So $\alpha (t,a,x)=a(e^{a}-1)\sim
a^{2}$ as $a\rightarrow 0.$ The drift coefficient does not contribute to the
non degeneracy condition (which is analogous to the uniform ellipticity
condition).

The paper is organized as follows. In Section 2 we give an integration by
parts formula of Malliavin type. This is analogous to the integration by
parts formulas given in \cite {BC} and \cite{BBM}. But there are two specific points:
first of all the integration by parts formula take into account the border
terms (in the above mentioned papers the border terms cancel because one
makes use  of some weights which are null on the border; but in the paper of
Kulik \cite{K2} such border terms appear as well). The second point is that we use
here a "one shot" integration by parts formula: in the classical gaussian Malliavin
calculus one employs all the noise which is
available - so one derives an infinite dimensional differential calculus
based on "all the increments" of the Brownian motion. The analogous approach
in the case of Poisson point measures is to use all the noise which comes
from the random structure (jumps). And this is the point of view of almost
all the papers on this topic. But in our paper we use just "one jump time"
which is chosen in a cleaver way (according to the non degeneracy
condition). In Section 3 we apply the general integration by parts formula
to stochastic equations with jumps. The basic noise is given by the jump
times.

\section{Integration by parts formula}
\subsection{Notations-derivative operators}

The abstract framework is quite similar  to the one developed  in \cite{BC} but we introduce here some modifications in order
to take into account the border terms appearing in the integration by parts formula. We consider a sequence of random variables $(V_{i})_{i\in \mathbb{N}^*}$ on a probability space $(\Omega ,\mathcal{F},P)$,
a sub $\sigma$-algebra $\mathcal{G}\subseteq \mathcal{F}$ and a random variable $J$, $\mathcal{G}$ measurable, with values in $\mathbb{N}$. 
Our aim is to establish a differential calculus based on the variables $(V_i)$, conditionally on $\mathcal{G}$. 
In order to derive an integration by parts formula, we need some  assumptions on the random variables $(V_i)$.
The main hypothesis  is that conditionally on $\mathcal{G},$  the law of $V_i$ admits a locally smooth density with respect to the Lebesgue measure.

\textbf{H0.} i) Conditionally on $\mathcal{G}$,  the random variables $(V_i)_{1 \leq i \leq J}$ are independent and for each $i \in \{1, \ldots, J\}$ the law of $V_i$
 is absolutely continuous with respect to the Lebesgue measure. We note $p_i$ the conditional density.

\hspace{1cm} ii) For all $i \in \{1, \ldots, J\}$, there exist some $\mathcal{G}$ measurable random variables $a_i$ and $b_i$  such that $-\infty<a_i<b_i<+ \infty$,  $(a_i,b_i)\subset \{p_i>0 \}$. We also assume  that $p_{i}$ admits a continuous  bounded derivative on $(a_i,b_i)$ and that  $\ln p_i$ is bounded on $(a_i, b_i)$.

We  define now the class of functions on which this differential calculus
will apply. We consider in this paper functions $f :\Omega \times \mathbb{R}^{\mathbb{N}^*}\rightarrow \mathbb{R}$  which can be written as 
\begin{equation}
f(\omega ,v)=\sum_{m=1}^{\infty }f^{m}(\omega
,v_{1},...,v_{m})1_{\{J(\omega )=m\}} \label{defS}
\end{equation}
where $f^{m}:\Omega \times \mathbb{R}^{m}\rightarrow \mathbb{R}$ are $\mathcal{G\times B}
(\mathbb{R}^{m})\mathcal{-}$measurable functions. 

In the following, we fix $L \in \mathbb{N}$ and we will perform integration by parts $L$ times. But we will use another set of variables for
each integration by parts. So for $1 \leq l \leq L$, we fix a set of indices $I_l \subset \{1, \ldots, J\}$ such that if $l \neq l'$, $I_l \cap I_{l'}=\emptyset$. In order to do $l$ integration by parts, we will use successively the variables $V_i, i \in I_l$ then the variables $V_i, i \in I_{l-1}$ and end with $V_i, i \in I_1$. Moreover, given $l$ we fix a partition $(\Lambda_{l,i})_{i \in I_l}$ of $\Omega$ such that the  sets $\Lambda_{l,i} \in \mathcal{G}, i \in I_l$.  If $\omega \in \Lambda_{l,i}$, we will use only the variable $V_i$  in our integration by parts.

With these notations, we define our basic spaces.  We consider in this paper random variables $F=f(\omega, V)$ where $V=(V_i)_i$ and $f$ is given by (\ref{defS}). To simplify the notation we write $F=f^J(\omega, V_1, \ldots ,V_J)$ so that conditionally on $\mathcal{G}$ we have $J=m$
and $F=f^m(\omega,V_1, \ldots, V_m)$.
We denote by $\mathcal{S}^0$  the space of random variables $F=f^J(\omega, V_1, \ldots, V_J)$ where $f^J$ is  a continuous function  on 
$O_J=\prod_{i=1}^J (a_i,b_i)$ such that there exists a $\mathcal{G}$ measurable random variable $C$ satisfying
\begin{equation}
\sup_{v \in O_J} \vert f^J( \omega, v) \vert  \leq C(\omega) < + \infty \quad {\mbox a.e.} \label{hypB}
\end{equation}
We also assume  that $f^J$ has left limits (respectively right  limits) in $a_{i}$
(respectively in $b_{i}$). Let us be more precise. 

With the notations $
V_{(i)}=(V_{1},...,V_{i-1},V_{i+1},...,V_{J})$  and $(V_{(i)}, v_i)=(V_1, \ldots, V_{i-1}, v_i, V_{i+1}, \ldots, V_J)$ for $v_i \in (a_i, b_i)$ our assumption
is that the following limits exist and are finite: 
\begin{equation}
\lim_{\varepsilon \rightarrow 0}f^{J}(\omega ,V_{(i)},a_i + \varepsilon):=F(a_{i}^+),\quad
\lim_{\varepsilon \rightarrow 0}f^{J}(\omega ,V_{(i)},b_i - \varepsilon):=F(b_{i}^-). \label{lim}
\end{equation}

Now for $k \geq 1$,
$\mathcal{S}^k(I_l)$ denotes the space of random variables  $F=f^J(\omega, V_1, \ldots, V_J) \in \mathcal{S}^0$, such that $f^J$ admits partial derivatives up to order $k$ with respect to the variables $v_i, i \in I_l$  and these partial derivatives belong to $\mathcal{S}^0$.

We are now able to define 
our differential operators.

$\square $ \textbf{The derivative operators.} We define  
$D_l:\mathcal{S}^1(I_l)\rightarrow 
\mathcal{S}^0(I_l):$ by
\begin{equation*}
D_l F:= 1_{O_J}(V)\sum_{i \in I_l} 1_{\Lambda_{l,i}}(\omega) \partial _{v_i}f(\omega ,V),      
\end{equation*}
where $O_J= \prod_{i=1}^J (a_i,b_i)$.

$\square $ \textbf{The divergence operators}  We note
\begin{equation}
p_{(l)}=\sum_{i \in I_l} 1_{\Lambda_{l,i}}p_i , \label{pl}
\end{equation}
and we define 
$\delta_l:\mathcal{S}^1(I_l)\rightarrow 
\mathcal{S}^0(I_l)$ 
by
\begin{eqnarray*}
\delta_l(F) = D_lF+F D_l \ln p_{(l)} = 1_{O_J}(V)\sum_{i \in I_l} 1_{\Lambda_{l, i}}(\partial _{v_i}F+F  \partial _{v_i} \ln p_i)
\end{eqnarray*}
We can easily see that if $F,U \in \mathcal{S}^1(I_l)$  we have
\begin{equation}
\delta_l(FU)=F\delta_l(U)+UD_lF .  \label{1.4}
\end{equation}

$\square $ \textbf{The border terms } \  Let  $U \in \mathcal{S}^0(I_l)$. We define (using the notation (\ref{lim}) )
\begin{eqnarray*}
[U]_l&=&\sum_{i \in I_l}1_{\Lambda_{l,i}} 1_{O_{J,i}}(V_{(i)})((U p_i)(b_i^-) -( U p_i)(a_i^+))
\end{eqnarray*}
with $O_{J,i}=\prod_{1 \leq j \leq J, j \neq i} (a_j,b_j)$

\subsection{Duality and basic integration by parts formula}

In our framework the duality between $\delta_l $ and $D_l$ is given by the
following proposition. In the sequel, we denote by 
$E_{\mathcal{G}}$ the conditional expectation with respect to the sigma-algebra $\mathcal{G}$.

\begin{prop}\label{duality}
Assuming  H0 then $\forall F, U\in 
\mathcal{S}^1(I_l)$  we have 
\begin{equation}
E_{\mathcal{G}}(UD_l F)=-E_{\mathcal{G}
}(F\delta_l (U)) + E_{\mathcal{G}} [FU]_l.  \label{1.2}
\end{equation}
\end{prop}
For simplicity, we assume in this proposition that the random variables $F$ and $U$ take value in $\mathbb{R}$ but such a result can
easily be extended to $\mathbb{R}^d$ value  random variables.
\begin{pf}
We have $E_{\mathcal{G}}(U D_lF)= \sum_{i \in I_l} 1_{\Lambda_{l,i}}E_{\mathcal{G}} 1_{O_J}(V)(\partial _{v_i}f^J(\omega ,V) u^J(\omega,V))$.
From H0  we obtain
\begin{equation*}
E_{\mathcal{G}} 1_{O_J}(V)( \partial _{v_i}f^J(\omega ,V) u^J(\omega,V))=E_{\mathcal{G}} 1_{O_{J,i}}(V_{(i)}) \int_{a_i}^{b_i} \partial
_{v_i}(f ^J) u^J p_{i}(v_i)dv_{i} . 
\end{equation*} 
By using the classical integration by parts formula, we have
$$
\int_{a_i}^{b_i} \partial
_{v_i}(f^J) u^J p_{i}(v_i)dv_i=[f^J u^J p_i]_{a_i}^{b_i}-\int_{a_i}^{b_i}f^J\partial_{v_i}(u^J p_i)dv_i.
$$
Observing that $\partial_{v_i}(u^J p_i)=( \partial_{v_i}(u^J ) +u^J    \partial_{v_i}( \ln p_i ))p_i$,
we have 
$$
E_{\mathcal{G}}(1_{O_J} (V)\partial _{v_i}f^J u^J)=E_{\mathcal{G}}1_{O_{J,i}}[(V_{(i)})f^J u^J p_i]_{a_i}^{b_i}-E_{\mathcal{G}} 1_{O_J}(V)F (\partial_{v_i}(U) +U \partial_{v_i}( \ln p_i ))
$$
and the proposition is proved. 

\end{pf}

We can now state a first integration by parts formula. 

\begin{prop} \label{IPP}
Let H0 hold true and  let $F \in 
\mathcal{S}^2(I_l)$, $G\in \mathcal{S}^1(I_l)$ and $\Phi :\mathbb{R}\rightarrow \mathbb{R}$ be a   bounded function with
bounded derivative. We assume that $F=f^J(\omega,V)$ satisfies the condition
\begin{equation}
\min_{ i \in I_l} \inf_{v \in O_J} \vert \partial_{v_i} f^J(\omega,v) \vert \geq \gamma(\omega) ,  \label{NDl}
\end{equation}
where $\gamma$ is $\mathcal{G}$ measurable and we define on $\{ \gamma>0 \}$
$$
(D_lF)^{-1}=1_{O_J}(V) \sum_{i \in I_l }1_{\Lambda_{l,i}} \frac{1}{\partial_{v_i} f(\omega,V)},
$$ 
then 
\begin{equation}
1_{\{ \gamma>0 \}} E_{\mathcal{G}}(\Phi^{(1)} (F)G)=-1_{\{ \gamma>0 \}}E_{\mathcal{G}}\left(\Phi
(F)H_l(F,G)\right) +  1_{\{ \gamma>0 \}}E_{\mathcal{G}}[\Phi(F)G (D_l F)^{-1}]_l \label{IPP1}
\end{equation}
with 
\begin{eqnarray}
H_l(F,G) &=&\delta_l (G (D_lF)^{-1})
=G \delta_l((D_lF)^{-1})+D_lG (D_lF)^{-1}.  \label{weight} 
\end{eqnarray}
\end{prop}

\begin{pf} We  observe that
$$
D_l\Phi (F)=1_{O_J}(V)\sum_{i \in I_l}1_{\Lambda_{l,i}} \partial_{v_i} \Phi(F)=1_{O_J}(V) \Phi^{(1)} (F)\sum_{i \in I_l}1_{\Lambda_{l,i}} \partial_{v_i} F,
$$
so that
\begin{eqnarray*}
D_l\Phi (F).D_lF&=&\Phi^{(1)}(F) (D_l F)^2,
\end{eqnarray*}
and then $1_{\{ \gamma>0 \}}\Phi^{(1)} (F) = 1_{\{ \gamma>0 \}} D_l \Phi(F). (D_l F)^{-1}$. Now since $F \in \mathcal{S}^2(I_l)$, we deduce that $(D_l F)^{-1} \in \mathcal{S}^1(I_l)$ on $\{ \gamma>0 \}$ and applying Proposition \ref{duality} with $U=G (D_l F)^{-1}$ we obtain the result.

\end{pf}

\subsection{Iterations of the integration by parts formula}

We will iterate the integration by parts formula given in Proposition \ref{IPP}. We recall that if we iterate $l$ times the integration by parts formula, we will  integrate by parts  successively with respect to the variables $(V_i)_{i \in I_k}$ for $1 \leq k \leq l$. In order to give some estimates of the weights appearing in these formulas
we introduce the following norm on $\mathcal{S}^l(\cup_{k=1}^l I_k)$, for $1 \leq l \leq L$.
\begin{equation}
\vert F \vert_{l}= \vert F \vert_{\infty} + \sum_{k=1}^l \sum_{1 \leq l_1<\ldots<l_k \leq l} \vert D_{l_1} \ldots D_{l_k} F \vert_{\infty}, \label{norml}
\end{equation}
where $\vert . \vert_{\infty}$ is defined on $\mathcal{S}^0$ by
$$
\vert F \vert_{\infty} = \sup_{v \in O_J} \vert f^J(\omega,v) \vert.
$$
For $l=0$, we set $\vert F \vert_0=\vert F \vert_{\infty} $.
We remark that  we have for $1 \leq l_1<\ldots<l_k \leq l$
$$
 \vert D_{l_1} \ldots D_{l_k} F \vert_{\infty} = \sum_{i_1 \in I_{l_1}, \ldots, i_k \in I_{l_k} } ( \prod_{j=1}^k1_{\Lambda_{l_j,i_j}}) \vert \partial_{v_{i_1}} \ldots \partial_{v_{i_k}} F \vert_{\infty},
$$
and since for each $l$ $(\Lambda_{l,i})_{i \in I_l}$ is a partition of $\Omega$,  for $\omega$ fixed, the preceding sum has only one term not equal to zero.
This family of norms satisfies for $F \in \mathcal{S}^{l+1}(\cup_{k=1}^{l+1} I_k)$ :
\begin{equation}
\vert F \vert_{l+1}=\vert D_{l+1} F \vert_l +\vert  F \vert_l  \quad \mbox{ so} \quad \vert D_{l+1} F \vert_l \leq \vert F \vert_{l+1}. \label{P1N}
\end{equation}
 Moreover it is easy to check that if $F, G \in \mathcal{S}^l(\cup_{k=1}^l I_k)$
\begin{equation}
\vert FG\vert_l \leq C_l \vert F \vert_l \vert G \vert_l, \label{prod}
\end{equation}
where $C_l$ is a constant depending on $l$. Finally for any function $\phi \in \mathcal{C}^l(\mathbb{R}, \mathbb{R})$ we have
\begin{equation}
\vert \phi(F) \vert_l \leq C_l \sum_{k=0}^l \vert \phi^{(k)}(F) \vert_{\infty} \vert F \vert_l^k  \leq C_l \max_{0 \leq k \leq l} 
\vert \phi^{(k)}(F) \vert_{\infty} (1+ \vert F \vert_l^l). \label{CR}
\end{equation}

With these notations we can iterate the integration by parts formula.
\begin{theo} \label{IPPE}
Let H0 hold true  and let $\Phi: \mathbb{R} \mapsto \mathbb{R}$ a bounded function with bounded derivatives up to order $L$. Let $F =f^J(w, V) \in \mathcal{S}^1(\cup_{l=1}^L I_l)$ such that
\begin{equation}
\inf_{i \in \{1, \ldots, J \}} \inf_{v \in O_J} \vert \partial_{v_i} f^J(\omega,v) \vert \geq  \gamma( \omega), \quad \gamma \in [0,1]  \quad \mathcal{G} {\mbox measurable}\label{ND}
\end{equation}
then  we have for $ l \in \{1, \ldots, L\}$, $G \in \mathcal{S}^l(\cup_{k=1}^l I_k)$ and  $F \in \mathcal{S}^{l+1}(\cup_{k=1}^l I_k)$
\begin{equation}
1_{\{ \gamma>0 \}}\vert E_{\mathcal{G}} \Phi^{(l)}(F)G \vert \leq C_l \vert \vert \Phi \vert \vert_{\infty}1_{\{ \gamma>0 \}}E_{\mathcal{G}} \left( \vert G \vert_l(1+ \vert p \vert_0 )^l
\Pi_l(F) \right) \label{IPPl}
\end{equation}
where $C_l$ is a constant depending on $l$, $\vert \vert \Phi \vert \vert_{\infty}= \sup_x \vert \Phi(x) \vert$, $\vert p\vert_0=\max_{l=1, \ldots,L} \vert p_{(l)} \vert_{\infty}$ and $\Pi_l(F)$ is defined  on $\{ \gamma>0 \}$ by
\begin{equation}
\Pi_l(F)=\prod_{k=1}^l
(1+ \vert( D_k F)^{-1} \vert_{k-1})(1+ \vert \delta_k((D_kF)^{-1} )\vert_{k-1}). \label{produit}
\end{equation}
Moreover  we have the bound 
\begin{equation}
\Pi_l(F) \leq C_l \frac{(1+ \vert \ln p \vert_1)^l}{\gamma^{l(l+2)} }\prod_{k=1}^l(1 + \vert F \vert_{k}^{k-1} + \vert D_k F \vert_{k}^{k-1})^2, \label{BPil}
\end{equation}
where $\vert \ln p\vert_1=\max_{i=1, \ldots,J} \vert (\ln p_i)' \vert_{\infty}$.
\end{theo}

\begin{pf}
We proceed by induction.
For $l=1$, we have  from Proposition \ref{IPP} since  $G \in \mathcal{S}^1( I_1)$ and  $F \in \mathcal{S}^{2}( I_1)$
$$
1_{\{ \gamma>0 \}}E_{\mathcal{G}}(\Phi^{(1)} (F)G)=-1_{\{ \gamma>0 \}}E_{\mathcal{G}}\left(\Phi
(F)H_1(F,G)\right) + 1_{\{ \gamma>0 \}} E_{\mathcal{G}}[\Phi(F)G (D_1 F)^{-1}]_1 .
$$
We have on $\{ \gamma>0 \}$
$$
\begin{array}{lll}
\vert H_1(F,G) \vert  & \leq &  \vert G \vert \vert \delta_1( (D_1F)^{-1}) \vert + \vert D_1 G \vert  \vert (D_1F)^{-1} \vert, \\
 & \leq & ( \vert G \vert_{\infty} +\vert D_1 G \vert _{\infty})(1+ \vert (D_1 F)^{-1} \vert_{\infty})(1+ \vert \delta_1((D_1F)^{-1} )\vert_{\infty}) , \\
  & = & \vert G \vert_1(1+ \vert (D_1 F)^{-1} \vert_{0})(1+ \vert \delta_1((D_1F)^{-1} )\vert_{0}).
\end{array}
$$
Turning to the border term $[\Phi(F)G (D_1 F)^{-1}]_1$, we check that
$$
\begin{array}{lll}
\vert [\Phi(F)G (D_1 F)^{-1}]_1 \vert & \leq & 2 \vert \vert \Phi \vert \vert_{\infty} \vert G \vert_{\infty} \sum_{i \in I_1} 1_{ \Lambda_{1,i}}  \vert \frac{1}{ \partial_{v_i} F} \vert_{\infty}  \sum_{i \in I_1} 1_{ \Lambda_{1,i}} \vert p_i \vert_{\infty}, \\
 & \leq & 2 \vert \vert \Phi \vert \vert_{\infty} \vert G \vert_{0} \vert (D_1F)^{-1} \vert_0 \vert p \vert_0.
\end{array}
$$
This proves the result for $l=1$.

Now assume that Theorem \ref{IPPE} is true for $l \geq 1$ and let us prove it for $l+1$. By assumption we have 
 $G \in \mathcal{S}^{l+1}(\cup_{k=1}^{l +1}I_k) \subset \mathcal{S}^{1}(I_{l+1})$ and  
 $F \in \mathcal{S}^{l+2}(\cup_{k=1}^{l+1} I_k) \subset \mathcal{S}^{2}(I_{l+1})$.
 Consequently we can apply Proposition \ref{IPP} on $I_{l+1}$. This gives
\begin{equation}
1_{\{ \gamma>0 \}}E_{\mathcal{G}}(\Phi^{(l+1)} (F)G)=-1_{\{ \gamma>0 \}}E_{\mathcal{G}}\left(\Phi^{(l)}
(F)H_{l+1}(F,G)\right) +  1_{\{ \gamma>0 \}}E_{\mathcal{G}}[\Phi^{(l)}(F)G (D_{l+1} F)^{-1}]_{l+1}, \label{recl}
\end{equation}
with
$$
H_{l+1}(F,G)=  G \delta_{l+1}((D_{l+1}F)^{-1})+D_{l+1}G (D_{l+1}F)^{-1}, 
$$
$$
[\Phi^{(l)}(F)G (D_{l+1} F)^{-1}]_{l+1} = \sum_{i \in I_{l+1}}1_{\Lambda_{l+1,i}} 1_{O_{J,i}}(V_{(i)})\left(( \Phi^{(l)}(F)G \frac{1}{\partial_{v_i} F} p_i)(b_i^-) -(\Phi^{(l)}(F)G \frac{1}{\partial_{v_i} F}  p_i)(a_i^+)\right) .
$$
We easily see  that $H_{l+1}(F,G) \in \mathcal{S}^l(\cup_{k=1}^l I_k)$ and so using the induction hypothesis we obtain
$$
1_{\{ \gamma>0 \}}\vert E_{\mathcal{G}}\Phi^{(l)}
(F)H_{l+1}(F,G) \vert \leq C_l \vert \vert \Phi \vert \vert_{\infty}1_{\{ \gamma>0 \}}E_{\mathcal{G}} \vert H_{l+1}(F,G) \vert_l(1 + \vert p \vert_0)^l \Pi_l(F),
$$
and we just have to bound $\vert H_{l+1}(F,G) \vert_l $ on $\{ \gamma>0 \}$. But using successively  (\ref{prod}) and (\ref{P1N})
$$
\begin{array}{lll}
\vert H_{l+1}(F,G) \vert_l  &  \leq  & C_l ( \vert G \vert_l \vert  \delta_{l+1}((D_{l+1}F)^{-1}) \vert_l + \vert D_{l+1}G \vert _l \vert  (D_{l+1}F)^{-1}) \vert_l, \\
 & \leq & C_l \vert G \vert_{l+1}(1+ \vert  (D_{l+1}F)^{-1}) \vert_l) (1+\vert  \delta_{l+1}((D_{l+1}F)^{-1}) \vert_l ).
\end{array}
$$
This finally gives
\begin{equation}
\vert E_{\mathcal{G}}\Phi^{(l)}
(F)H_{l+1}(F,G) \vert \leq C_l \vert \vert \Phi \vert \vert_{\infty}E_{\mathcal{G}} \vert G\vert_{l +1} (1+ \vert p \vert_0)^l\Pi_{l+1}(F). \label{Hl}
\end{equation}
So we just have to prove a similar inequality for $ E_{\mathcal{G}}[\Phi^{(l)}(F)G (D_{l+1} F)^{-1}]_{l+1}$. This reduces to consider
\begin{equation}
E_{\mathcal{G}}\sum_{i \in I_{l+1}}1_{\Lambda_{l+1,i}} 1_{O_{J,i}}(V_{(i)})( \Phi^{(l)}(F)G \frac{1}{\partial_{v_i} F} p_i)(b_i^-) 
= \sum_{i \in I_{l+1}}1_{\Lambda_{l+1,i}}  p_i(b_i^-) E_{\mathcal{G}} 1_{O_{J,i}}(V_{(i)})( \Phi^{(l)}(F)G \frac{1}{\partial_{v_i} F} )(b_i^-) \label{reduc}
\end{equation}
since the other term can be treated similarly. Consequently we just have to bound 
$$
\vert E_{\mathcal{G}} 1_{O_{J,i}}(V_{(i)})( \Phi^{(l)}(F)G \frac{1}{\partial_{v_i} F} )(b_i^-) \vert.
$$
Since all variables satisfy (\ref{hypB}), we obtain from Lebesgue Theorem,  using the notation (\ref{lim})
$$
E_{\mathcal{G}}1_{O_{J,i}}(V_{(i)})( \Phi^{(l)}(F)G \frac{1}{\partial_{v_i} F} )(b_i^-) = \lim_{\varepsilon \rightarrow 0 } E_{\mathcal{G}}1_{O_{J,i}}(V_{(i)})( \Phi^{(l)}(f^J(V_{(i)},b_i- \varepsilon))(g^J \frac{1}{\partial_{v_i} f^J} )(V_{(i)},b_i-\varepsilon).
$$
To shorten the notation we write simply $F(b_i - \varepsilon)=f^J(V_{(i)},b_i- \varepsilon)$.

Now one can prove that if $U \in \mathcal{S}^{l'}(\cup_{k=1}^{l+1}I_k)$ for $1 \leq l' \leq l$ then $\forall i \in I_{l+1}$, $U(b_i- \varepsilon) \in
\mathcal{S}^{l'}(\cup_{k=1}^{l}I_k)$ and $\vert U( b_i- \varepsilon)\vert_{l'} \leq \vert U \vert_{l'}$.
We deduce then that $\forall i \in I_{l+1}$
 $F(b_i- \varepsilon) \in \mathcal{S}^{l+1}(\cup_{k=1}^l I_k)$ and that $(G \frac{1}{\partial_{v_i} F} )(b_i- \varepsilon) \in \mathcal{S}^{l}(\cup_{k=1}^l I_k)$ and from induction hypothesis 
 $$
 \begin{array}{ll}
 \vert E_{\mathcal{G}}( \Phi^{(l)}(F(b_i- \varepsilon))1_{O_{J,i}}(G \frac{1}{\partial_{v_i} F} )(b_i-\varepsilon)\vert &\leq C_l  \vert \vert \Phi \vert \vert_{\infty} E_{\mathcal{G}}
 \vert G (b_i- \varepsilon) \vert_l \vert \frac{1}{\partial_{v_i} F(b_i- \varepsilon)} \vert_l (1 + \vert p \vert_0)^l \Pi_l(F(b_i- \varepsilon)), \\
& \leq  C_l  \vert \vert \Phi \vert \vert_{\infty}
E_{\mathcal{G}} \vert G  \vert_l \vert \frac{1}{\partial_{v_i} F} \vert_l (1 + \vert p \vert_0)^l \Pi_l(F).
\end{array}
$$
Putting this in (\ref{reduc}) we obtain
\begin{eqnarray}
\vert E_{\mathcal{G}}\sum_{i \in I_{l+1}}1_{\Lambda_{l+1,i}} 1_{O_{J,i}}( \Phi^{(l)}(F)G \frac{1}{\partial_{v_i} F} p_i)(b_i^-) \vert  &\leq& C_l
 \vert \vert \Phi \vert \vert_{\infty}
E_{\mathcal{G}}\vert G\vert_{l } (1+ \vert p \vert_0)^l\Pi_{l} (F)\sum_{i \in I_{l+1}}1_{\Lambda_{l+1,i}}  \vert p_i \vert_{\infty}\vert  \frac{1}{\partial_{v_i} F}\vert_l , \nonumber \\
 & \leq &C_l
 \vert \vert \Phi \vert \vert_{\infty}
E_{\mathcal{G}}\vert G\vert_{l } (1+ \vert p \vert_0)^{l+1}\Pi_{l}(F) \vert (D_{l+1}F)^{-1} \vert_l. \label{bord}
\end{eqnarray}
Finally  plugging (\ref{Hl}) and (\ref{bord}) in (\ref{recl})
$$
\begin{array}{lll}
\vert E_{\mathcal{G}}(\Phi^{(l+1)} (F)G) \vert & \leq & C_l
\vert \vert \Phi \vert \vert_{\infty} \left( E_{\mathcal{G}} \vert G\vert_{l +1} (1+ \vert p \vert_0)^l\Pi_{l+1}(F)
+ E_{\mathcal{G}}\vert G\vert_{l } (1+ \vert p \vert_0)^{l+1}\Pi_{l}(F) \vert (D_{l+1}F)^{-1} \vert_l \right), \\
 & \leq & C_l
\vert \vert \Phi \vert \vert_{\infty} E_{\mathcal{G}} \vert G \vert_{l+1}(1+ \vert p \vert_0 )^{l+1}
\Pi_{l+1}(F),
\end{array}
$$
and inequality (\ref{IPPl}) is proved for $l+1$. This achieves the first part of the proof of Theorem \ref{IPPE}.

It remains to prove (\ref{BPil}). We assume that $\omega \in \{ \gamma>0 \}$.

Let $1 \leq k \leq l$. 
We first notice that combining (\ref{P1N}) and (\ref{prod}) we obtain 
\begin{eqnarray*}
\left\vert \delta _{k}(F)\right\vert _{k-1} &\leq &\left\vert F\right\vert
_{k}(1+\left\vert D_{k}\ln p_{(k)}\right\vert _{\infty})  , 
\end{eqnarray*}
since $p_{(k)}$ only depends on the variables $V_i, i \in  I_{k}$.
 So we deduce the bound
\begin{eqnarray}
\left\vert \delta _{k}((D_k F)^{-1})\right\vert _{k-1} &\leq &\left\vert (D_k F)^{-1}\right\vert
_{k}(1+\left\vert \ln p\right\vert _{1})  \label{deltan}. 
\end{eqnarray}
Now we have 
\begin{equation*}
\vert (D_k F)^{-1} \vert_{k-1}=\sum_{i \in I_{k}} 1_{\Lambda_{k,i} }\vert \frac{1}{\partial_{v_i} F} \vert_{k-1}
\end{equation*}
From (\ref{CR}) with $\phi(x)=1/x$
$$
 \vert \frac{1}{\partial_{v_i} F} \vert_{k-1}\leq C_k \frac{(1 + \vert F \vert_k^{k-1})}{\gamma^{k}},
$$
and consequently
\begin{equation}
\vert (D_k F)^{-1} \vert_{k-1}\leq C_k \frac{(1 + \vert F \vert_k^{k-1})}{\gamma^{k}}. \label{DkF}
\end{equation}
Moreover we have using successively  (\ref{P1N}) and (\ref{DkF})
$$
\begin{array}{lll}
\vert (D_k F)^{-1} \vert_k & = & \vert (D_k F)^{-1} \vert_{k-1} +\vert D_k (D_k F)^{-1} \vert_{k-1} , \\
 & \leq & C_k \left( \frac{(1 + \vert F \vert_k^{k-1})}{\gamma^{k}} + \frac{(1 + \vert D_k F \vert_k^{k-1})}{\gamma^{k+1}} \right), \\
 & \leq & C_k  \frac{(1 +\vert F \vert_k^{k-1} +\vert D_k F \vert_k^{k-1})}{\gamma^{k+1}}.
\end{array}
$$
Putting this in (\ref{deltan})
\begin{equation}
\left\vert \delta _{k}((D_k F)^{-1})\right\vert _{k-1} \leq  C_k  \frac{(1 +\vert F \vert_k^{k-1} +\vert D_k F \vert_k^{k-1})}{\gamma^{k+1}}(1+\left\vert \ln p\right\vert _{1}) . \label{DDkF}
\end{equation}
Finally from (\ref{DkF}) and (\ref{DDkF}), we deduce
$$
\Pi_l(F) \leq C_l \frac{(1+\left\vert \ln p\right\vert _{1})^l }{\gamma^{l(l+2)} } \prod_{k=1}^l (1 +\vert F \vert_k^{k-1} +\vert D_k F \vert_k^{k-1})^2,
$$
and Theorem \ref{IPPE} is proved.
\end{pf}

\section{Stochastic equations with jumps}

\subsection{Notations and hypotheses}

We consider a Poisson point process $p$ with measurable state space $(E , \mathcal{B}(E))$. We refer to Ikeda and  Watanabe \cite{IK}
for the notation. We denote by $N$ the counting measure associated to $p$ so $N_t(A):= N((0,t) \times A)=\# \{s<t; p_s \in A\}$.
The intensity measure is $dt \times d\mu(a)$ where $\mu$ is a sigma-finite measure on $(E , \mathcal{B}(E))$ and we fix an
non decreasing sequence
$(E_n)$ of subsets of $E$ such that $E=\cup_n E_n$, $\mu(E_n) < \infty$  and $\mu(E_{n+1}) \leq \mu(E_n) + K$ for all $n$ and for a constant $K>0$.

We consider the one dimensional stochastic equation
\begin{equation}
X_t=x + \int_0^t \int_E c(s,a, X_{s^-})dN(s,a) + \int_0^tg(s, X_s) ds. \label{eq}
\end{equation}
Our aim is to give sufficient conditions on the coefficients $c$ and $g$ in order to prove that the law of $X_t$ is absolutely continuous
with respect to the Lebesgue measure and has a smooth density. We make the following assumptions on the coefficients $c$ and $g$.

{\bf H1.} We assume that the functions $c$ and $g$ are infinitely differentiable with respect to the variables $(t,x)$ and that there exist 
a bounded function $\overline{c}$ and  a constant $\overline{g}$,  such that
$$
\forall (t,a,x) \quad \vert c(t,a,x) \vert \leq  \overline{c}(a)(1+ \vert x \vert), \quad  \sup_{l+l' \geq 1}\vert \partial^{l'}_t \partial^l_x c(t,a,x)\vert \leq \overline{c}(a);
$$
$$
\forall (t,x)\quad \vert g(t,x) \vert \leq  \overline{g}(1+ \vert x \vert), \quad \sup_{ l+l' \geq 1}\vert \partial_t^{l'}\partial^l_x g(t,x)\vert \leq \overline{g};
$$
We assume moreover that $\int_E \overline{c}(a) d \mu(a) < \infty$.

Under H1, equation (\ref{eq}) admits a unique solution.

{ \bf H2.} We assume that there exists a measurable function $\hat{c}: E \mapsto \mathbb{R}_+$ such that $\int_E \hat{c}(a) d \mu(a) < \infty$ and 
$$
\forall (t,a,x)\quad \vert \partial_x c(t,a,x)(1+\partial_x c(t,a,x))^{-1} \vert \leq \hat{c}(a).
$$
To simplify the notation we take $\hat{c}=\overline{c}$. 
 Under H2, the tangent flow associated to (\ref{eq}) is invertible. At last we give a non-degeneracy condition wich will imply (\ref{ND}). We denote by $\alpha$ the function
 \begin{equation}
 \alpha(t,a,x)= g(t,x)-g(t,x+c(t,a,x))+(g \partial_xc+ \partial_t c)(t,a,x). \label{alpha}
 \end{equation}
{\bf H3.} We assume that there exists a measurable function $\underline{\alpha}: E \mapsto \mathbb{R}_+$ such that
$$
\forall (t,a,x) \quad \vert \alpha(t,a,x) \vert \geq \underline{\alpha}(a)>0,
$$
$$
\forall n \int_{E_n} \frac{1}{ \underline{\alpha}(a) }d \mu(a) < \infty \quad \mbox{ and} \quad\liminf_n \frac{1}{\mu(E_n)} \ln \left( \int_{E_n} \frac{1}{ \underline{\alpha}(a) }d \mu(a)\right) =\theta < \infty.
$$

We give in the following some examples where $E=(0,1]$ and $\underline{\alpha}(a)=a$.

\subsection{Main results and examples}

Following the methodology introduced in Bally and Cl\'ement \cite{BC}, our aim is to bound the Fourier transform of $X_t$, $\hat{p}_{X_t}( \xi )$, in terms of  $1/ \vert \xi \vert$, recalling that if $\int_{\mathbb{R} }\vert \xi \vert^q \vert \hat{p}_{X_t}( \xi) \vert d \xi < \infty$, for $q>0$,  then
the law of $X_t$ is absolutely continuous and its density is $\mathcal{C}^{[q]}$. This is done in the next proposition. The proof of this proposition relies on an approximation of $X_t$ which will be given in the next section.

\begin{prop} \label{fourier}
Assuming H1, H2 and H3 we have for all $n,L \in \mathbb{N}^*$
$$
\vert \hat{p}_{X_t}( \xi) \vert \leq C_{t,L} \left( e^{-\mu(E_n)t/(2L)}+ \frac{1}{\vert \xi\vert^L} A_{n,L} \right),
$$
with $A_{n,L}=\mu(E_n)^L (\int_{E_n} \frac{1}{\underline{\alpha}(a)} d \mu(a))^{L(L+2)}$.

\end{prop}

From this proposition, we deduce our main result.

\begin{theo} \label{density}
We assume that H1, H2 and H3 hold. Let $q \in \mathbb{N}$, then for $t> 16 \theta (q+2)(q+1)^2$, the law of $X_t$ is absolutely continuous with respect to the Lebesgue measure and its density is of class $\mathcal{C}^q$. In particular if $\theta=0$, the law of $X_t$ is absolutely continuous with respect to the Lebesgue measure and its density is of class $\mathcal{C}^{\infty}$ for every $t>0$.
\end{theo}

\begin{pf} From Proposition \ref{fourier}, we have
$$
\vert \hat{p}_{X_t}( \xi) \vert \leq C_{t,L} \left( e^{-\mu(E_n)t/2L}+ \frac{1}{\vert \xi\vert^L} A_{n,L} \right).
$$
Now  $\forall k,k_0>0$, if $t/2L >k \theta$, we deduce from H3 that for $n \geq n_L $ 
$$
t/2L >\frac{k}{\mu(E_n)} \ln (\int_{E_n} \frac{1}{\underline{\alpha}(a)} d \mu(a)) +\frac{k\ln \mu(E_n)}{k_0 \mu(E_n)}
$$
since the second term on the right hand side  tends to zero. This implies
$$
e^{\mu(E_n) t/2L} >(\int_{E_n} \frac{1}{\underline{\alpha}(a)} d \mu(a))^k \mu(E_n)^{k/k_0}.
$$ 
Choosing $k=L(L+2)$ and $k/k_0=L$, we obtain that for $n \geq n_L$  and $t/2L >L(L+2) \theta$
$$
e^{\mu(E_n) t/2L}>A_{n,L}.
$$
and then
$$
\begin{array}{lll}
\vert \hat{p}_{X_t}( \xi) \vert  & \leq  & C_{t,L} \left( e^{-\mu(E_n)t/2L}+\frac{1}{\vert \xi\vert^L} e^{\mu(E_n)t/2L} \right), \\
 & \leq & C_{t,L}(\frac{1}{B_n(t)}+\frac{B_n(t)}{\vert \xi \vert^L}),
\end{array}
$$
with $B_n(t)= e^{\mu(E_n)t/2L}$. Now recalling that $\mu(E_n)<\mu(E_{n+1} )\leq K+ \mu(E_n)$, we have $B_n(t)<B_{n+1}(t) \leq K_t B_n(t)$. Moreover since $B_n(t)$ goes to infinity with $n$ we have
$$
1_{\{ \vert \xi \vert^{L/2} \geq B_{n_L}(t) \}}= \sum_{n \geq n_L} 1_{ \{ B_n(t) \leq \vert \xi \vert^{L/2} < B_{n+1}(t)\}}.
$$ 
But if $B_n(t) \leq \vert \xi \vert^{L/2}< B_{n+1}(t)$, $\vert \hat{p}_{X_t}( \xi) \vert \leq C_{t,L}/ \vert \xi \vert^{L/2}$ and so
$$
\begin{array}{lll}
\int \vert \xi \vert^q  \vert \hat{p}_{X_t}( \xi) \vert d \xi  & = & \int_{\vert \xi \vert^{L/2}<B_{n_L(t)}}\vert \xi \vert^q  \vert \hat{p}_{X_t}( \xi)\vert d \xi
+  \int_{\vert \xi \vert^{L/2}\geq B_{n_L}(t)}\vert \xi \vert^q \vert \hat{p}_{X_t}( \xi) \vert d \xi, \\
 & \leq &C_{t,L, n_L} +  \int_{\vert \xi \vert^{L/2}\geq B_{n_L}(t)}\vert \xi \vert^{q-L/2}d \xi.
\end{array}
$$
For $q \in \mathbb{N}$, choosing $L$ such that $L/2-q>1$, we obtain  $\int \vert \xi \vert^q \vert \hat{p}_{X_t}( \xi) \vert d \xi < \infty$ for $t/2L >L(L+2) \theta$ and consequently the law of $X_t$ admits a density $\mathcal{C}^q$ for $t>2L^2(L+2) \theta$ and $L>2(q+1)$, that is $t>16\theta (q+1)^2(q+2)$ and Theorem \ref{density} is proved.

\end{pf}

We end this section with two examples

{\bf Example 1.} We take $E=(0,1]$, $\mu_{\lambda}=\sum_{k \geq 1} \frac{1}{k^{\lambda}} \delta_{1/k}$ with $0<\lambda<1$ and
$E_n=[1/n,1]$. We have $\cup_n E_n=E$, $\mu(E_n)= \sum_{k=1}^n \frac{1}{k^{\lambda}}$ and $\mu_{\lambda}(E_{n+1}) \leq \mu_{\lambda}(E_n) +1$ . We consider the process $(X_t)$ solution of (\ref{eq}) with $c(t,a,x)=a$ and $g(t,x)=g(x)$ assuming that the derivatives of $g$ are bounded
 and that $\vert g'(x) \vert \geq \underline{g}>0$. We have $\int_E a d \mu_{\lambda}(a)= \sum_{k \geq 1} \frac{1}{k^{\lambda+1}} < \infty$ so H1 and H2 hold.
 Moreover $\alpha(t,a,x)=g(x)-g(x+a)$ so $\underline{\alpha}(a)=\underline{g} a$. Now $\int_{E_n} \frac{1}{a} d \mu_{\lambda}(a) =\sum_{k=1}^n k^{1- \lambda}$ which is equivalent as $n$ go to infinity to $n^{2- \lambda}/(2- \lambda)$. Now we have
 $$
 \frac{1}{ \mu_{\lambda}(E_n)} \ln \left( \int_{E_n} \frac{1}{\underline{\alpha}(a)} d \mu_{\lambda}(a) \right)= \frac{ \ln( \underline{g} \sum_{k=1}^n k^{1- \lambda})}{ \sum_{k=1}^n \frac{1}{k^{\lambda}}} \thicksim_{ n \rightarrow \infty} C \frac{ \ln(n^{2- \lambda})}{n^{1- \lambda}} \rightarrow 0,
 $$
and then H3 is satisfied with $\theta=0$. We conclude from Theorem \ref{density} that $\forall t>0$, $X_t$ admits a density $\mathcal{C}^{\infty}$.

In the case $\lambda=1$, we have $\mu_1(E_n)=\sum_{k=1}^n \frac{1}{k} \thicksim_{n \rightarrow \infty} \ln n$  then
 $$
 \frac{1}{ \mu_{1}(E_n)} \ln \left( \int_{E_n} \frac{1}{\underline{\alpha}(a)} d \mu_{1}(a) \right)= \frac{ \ln( \underline{g} \sum_{k=1}^n 1)}{ \sum_{k=1}^n \frac{1}{k}} \thicksim_{ n \rightarrow \infty} 1,
 $$
and this gives H3 with $\theta=1$. So the density of $X_t$ is regular as soon as $t$ is large enough. In fact it is proved in Kulik \cite{K2} that
under some appropriate conditions the density of $X_t$ is not continuous for small $t$.

{\bf Example 2.} We take the intensity measure $\mu_{\lambda}$ as in the previous example and we consider the process $(X_t)$
solution of (\ref{eq}) with $g=1$ and $c(t,a,x)=ax$. This gives  $\overline{c}(a)=a$ and $\underline{\alpha}(a)=a$. So the conclusions are similar to example 1 in both cases $0< \lambda<1$ and $\lambda=1$.  But in this example we can compare our result to the one given by 
Ichikawa and Kunita \cite{IK}. They assume the condition
$$
\liminf_{u \rightarrow 0} \frac{1}{u^h} \int_{\vert a \vert \leq u} a^2 d \mu(a) >0, \quad (\star)
$$
for some $h \in (0,2)$. Here we have 
$$
\int_{\vert a \vert \leq u} a^2 d \mu(a) =\sum_{k \geq 1/u} \frac{1}{k^{2 + \lambda}}\thicksim_{u \rightarrow 0} \frac{u^{1+ \lambda}}{1 + \lambda}.
$$
So if $0<\lambda<1$, $(\star)$ holds and their results apply. In the case $\lambda=1$, $(\star)$ fails and they do not conclude. However
in our approach we conclude that the density of $X_t$ is $\mathcal{C}^q$ for $t>16(q+2)(q+1)^2$.

The next section is devoted to the proof of Proposition \ref{fourier}.

\subsection{Approximation of $X_t$ and integration by parts formula}
In order to bound the Fourier transform of the process $X_t$ solution of (\ref{eq}), we will apply the differential calculus developed in section 2. The first step consists in an approximation of $X_t$ by a random variable $X_t^N$ which can be viewed as an element of our
basic space $\mathcal{S}^0$. We assume that the process $(X_t^N)$ is solution of the discrete version of equation (\ref{eq})
\begin{equation}
X_t^N=x + \int_0^t \int_{E_N} c(s,a, X^N_{s^-})dN(s,a) + \int_0^tg(s, X^N_s) ds. \label{eqdis}
\end{equation}
Since $\mu(E_N) < \infty$, the number of jumps of the process $X^N$ on the interval $(0,t)$ is finite and consequently we may consider
the random variable $X_t^N$ as a function of these jump times and apply the methodology proposed in section 2. We denote by $(J_t^N)$ the Poisson process defined by $J_t^N=N((0,t), E_N)=\# \{s<t; p_s \in E_N \}$ and we note $(T_k^N)_{k \geq 1}$ its jump times. We also introduce the notation $\Delta_k^N=p_{T_k^N}$. With these notations, the process solution of (\ref{eqdis}) can be written 
\begin{equation}
X_t^N=x+\sum_{k=1}^{J_t^N} c(T_k^N, \Delta_k^N, X^N_{T_k^N-}) + \int_0^t g(s, X_s^N) ds . \label{eqdisbis}
\end{equation}
We will not work with all the variables $(T_k^N)_k$ but only with the jump times $(T_k^n)$ of the
Poisson process $J_t^n$, where $n<N$. In the following we will keep $n$ fixed and we will make $N$ go to infinity. We note $(T_k^{N,n})_k$ the jump times of the Poisson process $J_t^{N,n}=N((0,t), E_N\backslash E_n)$ and  $\Delta_k^{n,N}=p_{T_k^{n,N}}$. 
Now we fixe $L \in \mathbb{N}^*$, the number of integration by parts and we note $t_l=t l/L$, $0 \leq l \leq L$. Assuming that 
$J_{t_l}^n-J_{t_{l-1}}^n=m_l$ for $1 \leq l \leq L$, we denote by $(T_{l,i}^n)_{1 \leq i \leq m_l}$ the jump times of $J_t^n$ belonging to 
the time interval $(t_{l-1}, t_l)$.  In the following we assume that $m_l \geq 1$, $\forall l$. For $i=0$ we set $T_{l,0}^n=t_{l-1}$ and for $i=m_l+1$, $T_{l,m_l+1}^n=t_{l}$. With these definitions we choose our basic variables $(V_i, i \in I_l)$ as
\begin{equation}
(V_i, i \in I_l)=(T_{l,2i+1}^n, 0 \leq i \leq [(m_l-1)/2]). \label{Vi}
\end{equation}
The $\sigma$-algebra which contains the noise which is not involved in our differential calculus is 
\begin{equation}
\mathcal{G}=\sigma\{ (J_{t_l}^n)_{1 \leq l \leq L}; (T_{l,2i}^n)_{1 \leq2 i \leq m_l, 1 \leq l \leq L}; (T_k^{N,n})_k ; (\Delta_k^N)_k \}. \label{G}
\end{equation}
Using some well known results on Poisson processes, we easily see that conditionally on $\mathcal{G}$ the variables $(V_i)$ are independent and for $i \in I_l$ the law of $V_i$ conditionally on $\mathcal{G}$ is uniform on $(T_{l,2i}^n, T_{l,2i+2}^n)$ and we have
\begin{equation}
p_i(v)=\frac{1}{T_{l,2i+2}^n-T_{l,2i}^n} 1_{(T_{l,2i}^n, T_{l,2i+2}^n)}(v), \quad i \in I_l \label{pi},
\end{equation}
Consequently taking $a_i=T_{l,2i}^n$ and $b_i=T_{l,2i+2}^n$ we check that hypothesis H0 holds. It remains to define the localizing sets
$(\Lambda_{l,i})_{i \in I_l}$. 

We denote 
\[
h_{l}^{n}=\frac{t_{l}-t_{l-1}}{2m_{l}}=\frac{t}{2Lm_{l}}
\]
and $n_{l}=[(m_{l}-1)/2].$  We will work on the $\mathcal{G}$ measurable set 
\begin{equation}
\Lambda_{l}^n=\cup_{i=0}^{n_l} \{
T_{l,2i+2}^{n}-T_{l,2i}^{n}\geq h_{l}^{n}\}, \label{Ln}
\end{equation}
and we consider the following partition of this set:
\begin{eqnarray*}
\Lambda _{l,0} &=&\{T_{l,2}^{n}-T_{l,0}^{n}\geq h_{l}^{n}\}, \\
\Lambda _{l,i} &=&\cap
_{k=1}^{i}\{T_{l,2k}^{n}-T_{l,2k-2}^{n}<h_{l}^{n}\}\cap
\{T_{l,2i+2}^{n}-T_{l,2i}^{n}\geq h_{l}^{n}\},\quad i=1,...,n_{l}.
\end{eqnarray*}
After $L-l$ iterations of the integration by parts we will work with the variables $
V_{i},i\in I_{l}$ so the corresponding derivative is 
\[
D_{l}F=\sum_{i\in I_{l}}1_{\Lambda _{l,i}}\partial _{V_{i}}F=\sum_{i\in
I_{l}}1_{\Lambda _{l,i}}\partial _{T_{l,2i+1}^{n}}F.
\]
If we are on $\Lambda _{l}^n$ then we have at least one $i$ such that $
t_{l-1}\leq T_{l,2i}^{n}<T_{l,2i+1}^n<T_{l,2i+2}^{n}\leq t_{l}$ and $
T_{l,2i+2}^{n}-T_{l,2i}^{n}\geq h_{l}^{n}.$ Notice that in this case $
1_{\Lambda _{l,i}}\left\vert p_{i}\right\vert _{\infty }\leq (h_{l}^{n})^{-1}
$ and roughly speaking this means that the variable $V_{i}=T_{l,2i+1}^{n}$
gives a sufficiently large quantity of noise. Moreover, in order to perform $
L$ integrations by parts we will work on 
\begin{equation}
\Gamma _{L}^{n}=\cap _{l=1}^{L}\Lambda _{l}^n \label{GnL}
\end{equation}
and we will leave out the complementary of $\Gamma_{L}^{n}.$ The following
lemma says that on the set  $\Gamma _{L}^{n}$ we have enough noise
 and that the complementary of this set may be ignored. 

\begin{lem} \label{techni}
Using the notation given in Theorem \ref{IPPE} one has

i) $ \left\vert p\right\vert _{0} :=\max_{1\leq l\leq L}\sum_{i\in
I_{l}}1_{\Lambda _{l,i}}\left\vert p_{i}\right\vert _{\infty }\leq \frac{2L}{
t}  J_t^n$,

ii) $ P((\Gamma_{L}^{n})^{c}) \leq L\exp (-\mu (E_{n})t/2L)$.
\end{lem}

\begin{pf}
 As  mentioned before $1_{\Lambda _{l,i}}\left\vert
p_{i}\right\vert _{\infty }\leq (h_{l}^{n})^{-1}=2Lm_{l}/t\leq \frac{2L}{t
}J_{t}^{n}$ and so we have i). In order to prove ii) we have to
estimate $P((\Lambda_{l}^n)^{c})$ for $1 \leq l\leq L.$ We denote $s_{l}=\frac{1}{2}
(t_{l}+t_{l-1})$ and we will prove that $\{J_{t_{l}}^{n}-J_{s_{l}}^{n}\geq
1\}\subset \Lambda _{l}^n.$ Suppose first that $
m_{l}=J_{t_{l}}^{n}-J_{t_{l-1}}^{n}$ is even. Then $2n_{l}+2=m_{l}.$ If $
T_{l,2i+2}^{n}-T_{l,2i}^{n}< h_{l}^{n}$ for every $i=0,...,n_{l}$ then 
\[
T_{l,m_{l}}^n-t_{l-1}=\sum_{i=0}^{n_{l}}(T_{l,2i+2}^{n}-T_{l,2i}^{n})\leq
(n_{l}+1)\times \frac{t}{2Lm_{l}}\leq \frac{t}{4L}\leq s_{l}-t_{l-1}
\]
so there are no jumps in $(s_{l},t_{l}).$ Suppose now that $m_{l}$ is odd
so $2n_{l}+2=m_{l}+1$ and $T_{l,2n_{l}+2}^{n}=t_{l}.$ If we have $T_{l,2i+2}^{n}-T_{l,2i}^{n}< h_{l}^{n}$ for every $
i=0,...,n_{l},$ then we deduce
\[
\sum_{i=0}^{n_{l}}(T_{l,2i+2}^{n}-T_{l,2i}^{n})< (n_{l}+1)\times \frac{t}{
2Lm_{l}}< \frac{m_l+1}{m_l} \frac{t}{4L}\leq  \frac{t}{2L},
\]
and there are no jumps in $(s_{l},t_{l}).$
So we have proved that $\{J_{t_{l}}^{n}-J_{s_{l}}^{n}\geq 1\}\subset \Lambda
_{l}^n$ and since $P(J_{t_{l}}^{n}-J_{s_{l}}^{n}=0)=\exp (-\mu (E_{n})t/2L)$
the inequality $ii)$ follows. 
\end{pf}

Now we will apply Theorem \ref{IPPE}, with $F^N=X_t^N$, $G=1$ and $\Phi_{\xi}(x)=e^{i \xi x}$. So we have to check that
$F^N \in \mathcal{S}^{L+1}( \cup_{l=1}^L I_l)$ and that condition (\ref{ND}) holds. Moreover we have to bound $\vert F^N \vert_l^{l-1}$ 
and $\vert D_l F^N \vert_l^{l-1}$, for $1 \leq l \leq L$. This needs some preliminary lemma.

\begin{lem} \label{derdet}
Let $v=(v_i)_{i \geq 0}$ a positive non increasing sequence with $v_0=0$ and $(a_i)_{i \geq 1}$ a sequence of $E$. We define
$J_t(v)$  by $J_t(v)=v_i$ if $v_i \leq t <v_{i+1}$ and we consider the process solution of
\begin{equation}
X_t=x+ \sum_{k=1}^{J_t} c(v_k, a_k, X_{v_k-})+ \int_0^t g(s,X_s) ds \label{deteq}.
\end{equation}
We assume that H1 holds. Then $X_t$ admits some derivatives with respect to $v_i$ and if we note $U_i(t)= \partial_{v_i} X_t$ and 
$W_i(t)=\partial^2_{v_i} X_t$, the processes $(U_i(t))_{ t \geq v_i}$ and $(W_i(t))_{t \geq v_i}$ solve respectively
\begin{equation}
U_i(t)= \alpha(v_i, a_i, X_{v_i-}) + \sum_{k=i+1}^{J_t} \partial_x c(v_k, a_k, X_{v_k-})U_i(v_k-) + \int_{v_i}^t \partial_x g(s, X_s) U_i(s) ds
\label{DFN},
\end{equation}
\begin{equation}
W_i(t)= \beta_i(t) + \sum_{k=i+1}^{J_t} \partial_x c(v_k, a_k, X_{v_k-})W_i(v_k-) + \int_{v_i}^t \partial_x g(s, X_s) W_i(s) ds
\label{D2FN},
\end{equation}
with
$$
\begin{array}{cll}
\alpha(t,a,x) & = & g(t,x)-g(t,x+c(t,a,x))+ g(t,x) \partial_x c(t,a,x) + \partial_t c(t,a,x) , \\
\beta_i(t) & = & \partial_t \alpha(v_i,a_i ,X_{v_i-})+\partial_x \alpha(v_i,a_i ,X_{v_i-}) g(v_i, X_{v_i-})-\partial_x g(v_i, X_{v_i}) U_i(v_i) \\
 & &+ \sum_{k=i+1}^{J_t}\partial_x^2 c(v_k, a_k, X_{v_k-})(U_i(v_k-) )^2 +\int_{v_i}^t \partial_x^2 g(s, X_s) (U_i(s))^2 ds.
\end{array}
$$
\end{lem}

\begin{pf}
If $s<v_i$, we have $\partial_{v_i} X_s=0$. Now we have
$$
X_{v_i-}=x+\sum_{k=1}^{v_{i-1}} c(v_k, a_k ,X_{v_k-})+\int_0^{v_i }g(s,X_s) ds,
$$
and consequently
$$
\partial_{v_i}X_{v_i-}=g(v_i, X_{v_i-}).
$$
For $t>v_i$, we observe that
$$
X_t=X_{v_i-}+ \sum_{k=v_i}^{J_t} c(v_k, a_k ,X_{v_k-})+  \int_{v_i}^t g(s, X_s) ds,
$$
this gives
$$
\begin{array}{lll}
\partial_{v_i} X_t& = & g(v_i, X_{v_i-})+g(v_i, X_{v_i-}) \partial_x c(v_i, a_i,X_{v_i-}) + \partial_t c(v_i, a_i,X_{v_i-})-g(v_i, X_{v_i}) \\
 & & + \sum_{k=i+1}^{J_t} \partial_x c(v_k, a_k, X_{v_k-})\partial_{v_i} X_{v_k-} + \int_{v_i}^t \partial_x g(s, X_s) \partial_{v_i} X_s  ds.
\end{array}
$$
Remarking that $X_{v_i}=X_{v_i-}+c(v_i, a_i,X_{v_i-})$, we obtain (\ref{DFN}). The proof of (\ref{D2FN}) is similar and we omit it.

\end{pf}

We give next a bound for $X_t$ and its derivatives  with respect to the variables $(v_i)$.
\begin{lem} \label{Bdet}
Let $(X_t)$ the process solution of (\ref{deteq}).  We assume that H1 holds and we note 
$$
n_t(\overline{c})=\sum_{k=1}^{J_t} \overline{c}(a_k).
$$
Then we have:
$$
\sup_{s \leq t} \vert X_t\vert \leq C_t(1+n_t(\overline{c}))e^{n_t(\overline{c})}.
$$
Moreover
$\forall l \geq 1$, there exist some constants $C_{t,l}$ and $C_l$ such that
$\forall (v_{k_i})_{i=1, \ldots, l}$ with $t>v_{k_l}$, we have 
$$
\sup_{v_{k_l} \leq s \leq t} \vert \partial_{v_{k_1}} \ldots \partial_{v_{k_{l-1}}} U_{k_l}(s) \vert
+\sup_{v_{k_l} \leq s \leq t} \vert \partial_{v_{k_1}} \ldots \partial_{v_{k_{l-1}}} W_{k_l}(s) \vert \leq C_{t,l}(1+n_t(\overline{c}))^{C_l}e^{C_l n_t(\overline{c})}.
$$
\end{lem}
We observe that the previous bound does not depend on the variables $(v_i)$.
\begin{pf}
We just give a sketch of the proof.
We first remark that the process  $(e_t)$ solution of
$$
e_t=1+\sum_{k=1}^{J_t} \overline{c}(a_k) e_{v_k-} + \overline{g}\int_0^t e_s ds,
$$
is given by $e_t=\prod_{k=1}^{J_t} (1+\overline{c}(a_k)) e^{\overline{g}t}$.
Now from H1, we deduce for $s \leq t$
$$
\begin{array}{lll}
\vert X_s \vert & \leq & \vert x \vert + \sum_{k=1}^{J_s} \overline{c}(a_k)(1 + \vert X_{v_k-} \vert ) + \int_0^s \overline{g} (1 + \vert X_u \vert )du, \\
 & \leq &  \vert x \vert + \sum_{k=1}^{J_t} \overline{c}(a_k) +\overline{g} t + \sum_{k=1}^{J_s} \overline{c}(a_k) \vert X_{v_k-} \vert  + \int_0^s \overline{g}  \vert X_u \vert du, \\
  & \leq & ( \vert x \vert + \sum_{k=1}^{J_t} \overline{c}(a_k) +\overline{g} t ) e_s 
 \end{array}
$$
where the last inequality follows from Gronwall lemma.
Then using the previous remark
\begin{equation}
\sup_{s \leq t} \vert X_s\vert \leq C_t (1+ n_t(\overline{c}))\prod_{k=1}^{J_t} (1+ \overline{c}(a_k)) \leq C_t (1+ n_t(\overline{c})) e^{n_t(\overline{c})}. \label{BX}
\end{equation}
We check easily that $\vert \alpha(t,a,x) \vert \leq C(1+ \vert x \vert ) \overline{c}(a)$, and we get successively from (\ref{DFN}) and (\ref{BX})
$$
\sup_{v_{k_l} \leq s \leq t} \vert U_{k_l}(s) \vert  \leq C_t(1+ \vert X_{v_{k_l}-} \vert ) \overline{c}(a_{k_l})(1+ n_t(\overline{c})) e^{n_t(\overline{c})} \leq C_t (1+ n_t(\overline{c}))^2 e^{2 n_t(\overline{c})}.
$$
Putting this in (\ref{D2FN}), we obtain a similar bound for $\sup_{v_{k_l} \leq s \leq t} \vert W_{k_l}(s) \vert$ and we end the proof
of Lemma \ref{Bdet} by induction since we can derive equations for the higher order derivatives of $U_{k_l}(s)$ and $W_{k_l}(s)$
analogous to (\ref{D2FN}).

\end{pf}
 
 We come back to the process $(X_t^N)$ solution of (\ref{eqdis}). We recall that $F^N=X_t^N$ and we will check that $F^N$ satisfies
 the hypotheses of Theorem \ref{IPPE}.

\begin{lem} \label{verhyp}
i) We assume that H1 holds. Then $\forall l \geq 1$, $\exists C_{t,l}, C_l$ independent of $N$ such that
$$
\vert F^N \vert_l + \vert D_l F^N \vert_l \leq C_{t,l} \left( (1+N_t(\overline{c}))e^{N_t(\overline{c}) }\right)^{C_l},
$$
with  $N_t(\overline{c})= \int_0^t \int_E \overline{c}(a) dN(s,a) $.

ii) Moreover if we assume in addition that H2 and H3 hold and that $m_l=J_{t_l}^n-J_{t_{l-1}}^n \geq 1$, $\forall l \in \{1, \ldots, L\}$  then we have $\forall 1 \leq l \leq L$, $\forall i \in I_l$ 
$$
\vert \partial_{V_i}  F^N \vert \geq \left(e^{2 N_t(\overline{c})} N_t (1_{E_n} 1/\underline{\alpha} ) \right)^{-1} :=\gamma_n
$$
and  (\ref{ND}) holds. 

\end{lem}
We remark that on the non degeneracy set $\Gamma_L^n$ given by (\ref{GnL}) we have at least one jump on $(t_{l-1}, t_l)$, that is $m_l \geq 1$, $\forall l \in \{1, \ldots, L\}$. Moreover we have $\Gamma_L^n \subset \{\gamma_n>0 \}$.
\begin{pf}
The proof of i) is a straightforward consequence of Lemma \ref{Bdet}, replacing $n_t(\overline{c})$ by $\sum_{p=1}^{J_t^N} \overline{c}(\Delta_p^N)$ and observing that 
$$
\sum_{p=1}^{J_t^N} \overline{c}(\Delta_p^N)=  \int_0^t \int_{E_N} \overline{c}(a) dN(s,a) \leq \int_0^t \int_E \overline{c}(a) dN(s,a)=N_t(\overline{c}).
$$
Turning to ii) we have from Lemma \ref{derdet} 
$$
\partial_{T_k^N} X_t^N=\alpha(T_k^N, \Delta_k^N, X^N_{T_k^N-}) + \sum_{p=k+1}^{J_t^N} \partial_x c(T_p^N, \Delta_p^N, X^N_{T_p^N-}) \partial_{T_k^N} X_{T_p^N-}^N+ \int_{T_k^N}^t \partial_x g(s, X^N_s)  \partial_{T_k^N} X_s ds.
$$
Assuming H2, we define $(Y_t^N)_t$ and $(Z_t^N)_t$ as the solutions of the equations
$$
\begin{array}{lll}
Y_t^N & = & 1 + \sum_{p=1}^{J_t^N} \partial_x c(T_p^N, \Delta_p^N, X^N_{T_p^N-}) Y_{T_k^N-} + \int_0^t \partial_x g(s, X_s^N) Y_s^N ds, 
\\
Z_t^N & = & 1 - \sum_{p=1}^{J_t^N} \frac{\partial_x c(T_p^N, \Delta_p^N, X^N_{T_p^N-})}{1+\partial_x c(T_p^N, \Delta_p^N, X^N_{T_p^N-})} Z_{T_k^N-} -\int_0^t \partial_x g(s, X_s^N) Z_s^N ds.
\end{array}
$$
We have $Y_t^N \times Z_t^N=1$, $\forall t \geq 0$ and 
$$
\vert Y_t^N \vert \leq e^{t \overline{g}} e^{N_t(1_{E_N} \overline{c})} \leq e^{N_t( \overline{c})} ,\quad \vert Z_t^N \vert = \vert \frac{1}{ Y_t^N} \vert \leq e^{N_t( \overline{c})}.
$$
Now one can easily check that
$$
\partial_{T_k^N} X_t^N=\alpha(T_k^N, \Delta_k^N, X^N_{T_k^N-}) Y_t^N Z^N_{T_k^N},
$$
and using H3 and the preceding bound it yields
$$
\vert \partial_{T_k^N} X_t^N\vert \geq  e^{-2N_t( \overline{c})} \underline{\alpha}(\Delta_k^N).
$$
Recalling that we do not consider the derivatives with respect to all the variables $(T_k^N)$ but only with respect to $(V_i)=(T^n_{l,2i+1})_{l,i}$
with $n<N$ fixed, we have $\forall 1 \leq l \leq L$ and $ \forall i \in I_l$
$$
\vert \partial_{V_i} X_t^N\vert \geq e^{-2N_t( \overline{c})} \left(\sum_{p=1}^{J_t^n} \frac{1}{\underline{\alpha}(\Delta_p^n)}\right)^{-1}=
\left(e^{2 N_t(\overline{c})} N_t (1_{E_n} 1/\underline{\alpha} ) \right)^{-1},
$$
and Lemma \ref{verhyp} is proved.

\end{pf}

With this lemma we are at last able to prove Proposition \ref{fourier}.

\noindent
{\bf Proof of Proposition \ref{fourier}: }
From Theorem \ref{IPPE} we have since $\Gamma_L^n \subset \{ \gamma_n >0 \}$
$$
1_{\Gamma_L^n} \vert E_{\mathcal{G} }\Phi^{(L)}(F^N) \vert \leq C_L \vert \vert \Phi \vert \vert_{\infty} 1_{\Gamma_L^n} E_{\mathcal{G}} (1+\vert p_0 \vert)^L \Pi_L(F^N).
$$
Now from Lemma \ref{techni} i) we have 
$$
\vert p_0 \vert  \leq 2 L J_t^n/t
$$
and moreover we can check that $\vert \ln p \vert_1=0$. So we deduce from Lemma \ref{verhyp}
$$
\Pi_L(F^N) \leq \frac{C_{t,L}}{\gamma_n^{L(L+2)} }\left( (1+N_t(\overline{c}))e^{N_t(\overline{c}) }\right)^{C_L} \leq 
C_{t,L}N_t (1_{E_n} 1/\underline{\alpha} )^{L(L+2)} \left( (1+N_t(\overline{c}))e^{N_t(\overline{c}) }\right)^{C_L} .
$$
This finally gives
\begin{equation}
\vert E1_{\Gamma_L^n} \Phi^{(L)}(F^N) \vert \leq \vert \vert \Phi \vert \vert_{\infty} C_{t,L} E\left((J_t^N)^L N_t (1_{E_n} 1/\underline{\alpha} )^{L(L+2)} \left( (1+N_t(\overline{c}))e^{N_t(\overline{c}) }\right)^{C_L} \right). \label{Bint}
\end{equation}
Now we know from a classical computation (see for example \cite{BC}) that the Laplace transform of $N_t(f)$
satisfies
\begin{equation}
E e^{-s N_t(f)}=e^{-t \alpha_{f}(s)}, \quad \alpha_{f}(s)= \int_E (1-e^{-s f(a)}) d \mu(a). \label{laplace}
\end{equation}
From H1,  we have $\int_E \overline{c}(a) d \mu(a)< \infty$, so we deduce using (\ref{laplace}) with $f= \overline{c}$ that, $\forall q>0$
$$
E \left( (1+N_t(\overline{c}))e^{N_t(\overline{c}) }\right)^q \leq C_{t,q} < \infty.
$$
Since $J_t^n$ is a Poisson process with intensity $t \mu(E_n)$, we have $\forall q>0$
$$
E (J_t^n)^q \leq C_{t,q} \mu(E_n)^q.
$$
Finally, using once again (\ref{laplace}) with $f=1_{E_n} 1/\underline{\alpha}$ we see easily that $\forall q>0$
$$
EN_t (1_{E_n} 1/\underline{\alpha} )^q \leq C_{t,q} \left( \int_{E_n} \frac{1}{\underline{\alpha}(a)} d \mu(a) \right)^q.
$$
Turning back to (\ref{Bint}) and combining Cauchy-Schwarz inequality  and the previous bounds we deduce
\begin{equation}
\vert E1_{\Gamma_L^n} \Phi^{(L)}(F^N) \vert \leq \vert \vert \Phi \vert \vert_{\infty} C_{t,L} \mu(E_n)^L \left(  \int_{E_n} \frac{1}{\underline{\alpha}(a)} d \mu(a)\right)^{L(L+2)} = \vert \vert \Phi \vert \vert_{\infty}C_{t,L} A_{n,L}. \label{Bfin}
\end{equation}
We are now ready to give a bound for $\hat{p}_{X_t^N}(\xi)$.
We have $\hat{p}_{X_t^N}(\xi)= E \Phi_{\xi}(F^N)$, with $\Phi_{\xi}(x)=e^{i \xi x}$. Since $\Phi^{(L)}_{\xi}(x)=(i \xi)^L \Phi_{\xi}(x)$,
we can write $\vert \hat{p}_{X_t^N}(\xi)\vert=\vert E\Phi^{(L)}_{\xi}(F^N) \vert / \vert \xi \vert^L$ and consequently we deduce from (\ref{Bfin})
$$
\vert \hat{p}_{X_t^N}(\xi)\vert \leq P((\Gamma_L^n)^c) +C_{t,L} A_{n,L}/ \vert \xi \vert^L.
$$
But from Lemma \ref{techni} ii)
we have
$$
P((\Gamma_L^n)^c)  \leq L e^{-\mu(E_n) t/(2L)}
$$
and finally
$$
\vert \hat{p}_{X_t^N}(\xi) \vert  \leq C_{L,t} \left( e^{-\mu(E_n) t/(2L)} + A_{n,L}/ \vert \xi \vert^L \right).
$$
We achieve the proof of  Proposition \ref{fourier} by letting $N$ go to infinity, keeping $n$ fixed. 

{\hfill\mbox{$\diamond$}\medskip}

\end{document}